\documentclass{ifacconf}

\usepackage{graphicx}      
\usepackage{natbib}        
\usepackage{cases}
\usepackage{amsmath}
\usepackage{amsfonts}
\usepackage{amssymb}
\usepackage{latexsym}
\usepackage{xcolor}

\DeclareMathOperator{\Tr}{Tr}
\newcommand{\bmat}[1]{\begin{bmatrix} #1 \end{bmatrix}}

\begin{document}

\begin{frontmatter}

\title{Adaptive Optimal PMU Placement Based on Empirical Observability Gramian\thanksref{footnoteinfo}}

\thanks[footnoteinfo]{This work was supported in part by U.S. Department of Energy, Office of Electricity Delivery and Energy Reliability under contract DE-AC02-06CH11357, 
the CURENT engineering research center, and Naval Research Laboratory and Defense Advanced Research Projects Agency.}

\author[First]{Junjian Qi} 
\author[Second]{Kai Sun} 
\author[Third]{Wei Kang}

\address[First]{Energy Systems Division, Argonne National Laboratory, Argonne, IL 60439 USA (e-mail: jqi@anl.gov)}
\address[Second]{Department of EECS, University of Tennesee, Knoxville, TN 37996 USA (e-mail: kaisun@utk.edu)}
\address[Third]{Department of Applied Mathematics, Navel Postgraduate School, Monterey, CA 93943 USA (e-mail: wkang@nps.edu)}

\begin{abstract}                
In this paper, we compare four measures of the empirical observability gramian, including the determinant, the trace, the minimum eigenvalue, and the condition number, which can be used to quantify the observability of system states and to obtain the optimal PMU placement for power system dynamic state estimation. An adaptive optimal PMU placement method is proposed by automatically choosing proper measures as the objective function. It is shown that when the number of PMUs is small and thus the observability is very weak, the minimum eigenvalue and the condition number are better measures of the observability and are preferred to be chosen as the objective function. The effectiveness of the proposed method is validated by performing dynamic state estimation on an Northeast Power Coordinating Council (NPCC) 48-machine 140-bus system with the square-root unscented Kalman filter.
\end{abstract}


\end{frontmatter}

\section{Introduction}

Increasing integration of intermittent renewable energy and current effort of developing smart grid will make electric power systems more and more dynamic. However, the most widely studied power system static state estimation (SSE) (\cite{Schweppe,Abur,Monticelli,Irving,He,Qi}) cannot capture the dynamics of power systems well due to its dependency on slow update rates of Supervisory Control and Data Acquisition (SCADA) systems.

By contrast, real-time dynamic state estimation (DSE) enabled by phasor measurement units (PMUs), which has high update rates and 
high global positioning system (GPS) synchronization accuracy, can provide accurate dynamic states of the system and 
thus will play a critical role in achieving real-time wide-area monitoring, protection, and control (\cite{Begovic,Qi3}).
Until now DSE has been implemented by extended Kalman filter (\cite{Huang,Ghahremani}), 
unscented Kalman filter (\cite{Wang,Singh}), square-root unscented Kalman filter (\cite{Qi2,Qi1}), extended particle filter (\cite{Zhou}), 
cubature Kalman filter(\cite{Taha}), and observers (\cite{Taha1,Taha}).

The well-known optimal PMU placement (OPP) problem was originally developed for SSE. 
It is mainly based on the topological observability criterion, which only specifies that the power system states should be uniquely estimated with the minimum number of PMUS but neglects important parameters such as transmission line admittances by only focusing on the binary connectivity graph (\cite{Baldwin,Li}). Under this framework, many approaches have been proposed, such as mixed integer programming (\cite{Xu,Gou}), binary search (\cite{Chakrabarti}), metaheuristics (\cite{Milosevic,Aminifar}), particle swarm optimization (\cite{Chakrabarti1}), and eigenvalue-eigenvector based approaches (\cite{Almutairi,Korba}). 
An information-theoretic criterion is also proposed to generate highly informative PMU configurations (\cite{Li}).

However, not much research has been done on OPP for DSE. In (\cite{Kamwa}) numerical PMU configuration algorithms are proposed to maximize the overall sensor response
while minimizing the correlation among sensor outputs based on the system response under many contingencies. 
In (\cite{Sun}) an OPP strategy is proposed to ensure a satisfactory state tracking performance but it depends on a specific Kalman filter. In (\cite{Qi1}) the empirical observability gramian (\cite{Lall,Lall1,Hahn,Singh1}) is applied to quantify the degree of observability of the system states and OPP is achieved by maximizing the determinant of the empirical observability gramian. Compared with (\cite{Kamwa}) and (\cite{Sun}), (\cite{Qi1}) has a quantitative measure of observability, which makes it possible to optimize PMU locations from the point view of the observability of nonlinear systems. Besides, it only needs to deal with the system under typical power flow conditions and also does not depend on the specific realization of any Kalman filter.

As in (\cite{Qi1}), there are various measures of the empirical observability gramian that can be used to quantify the observability of the system state. 
In this paper, we compare these measures and further propose an adaptive OPP method for power system dynamic state estimation by automatically choosing proper measures as the objective function to guarantee best observability.

The remainder of this paper is organized as follows. 
Section \ref{estimation} briefly introduces power system dynamic state estimation. 
Section \ref{gramian} discusses the definition and implementation of empirical observability gramian. 
Section \ref{adaptive} introduces the OPP method based on the empirical observability gramian and proposes an adaptive OPP method by making full use of different measures. The results are presented in Section \ref{case} in order to test and validate the proposed method. 
Finally the conclusion is drawn in Section \ref{conclusion}.


\section{Power System Dynamic State Estimation} \label{estimation}

Different from SSE, DSE estimates the dynamic states (internal states of generators), rather than the static states (voltage magnitude and phase angles of buses). 
In order to perform DSE, the nonlinear dynamics and the outputs of a power system is described in the following form:
\begin{equation}~\label{n1}
\left\{
\arraycolsep=1.4pt\def\arraystretch{1.2}
\begin{array}{ll}
\dot{\boldsymbol{x}}=\boldsymbol{f}(\boldsymbol{x},\boldsymbol{u}) \\
\boldsymbol{y}=\boldsymbol{h}(\boldsymbol{x},\boldsymbol{u})
\end{array}
\right.
\end{equation}
where $\boldsymbol{f}(\cdot)$ and $\boldsymbol{h}(\cdot)$ are the state transition and output functions, $\boldsymbol{x} \in \mathbb{R}^n$ is the state vector, $\boldsymbol{u} \in \mathbb{R}^v$ is the input vector, and $\boldsymbol{y}\in \mathbb{R}^p$ is the output vector.

We consider two types of generator models, the fourth-order transient model and second-order classical model. 
For generator with transient model, the fast sub-transient dynamics and saturation effects are ignored and the generator is described 
by fourth-order differential equations: 
\begin{equation}~\label{gen model}
\left\{
\arraycolsep=1.4pt\def\arraystretch{2.0}
\begin{array}{ll}
\dot{\delta_i}=\omega_i-\omega_0 \\
\dot{\omega}_i=\frac{\omega_0}{2H_i}\Big(T_{mi}-T_{ei}-\frac{K_{Di}}{\omega_0}(\omega_i-\omega_0)\Big) \\
\dot{e}'_{qi}=\frac{1}{T'_{d0i}}\Big(E_{fdi}-e'_{qi}-(x_{di}-x'_{di})\,i_{di}\Big) \\
\dot{e}'_{di}=\frac{1}{T'_{q0i}}\Big(-e'_{di}+(x_{qi}-x'_{qi})\,i_{qi}\Big)
\end{array}
\right.
\end{equation}
where $i$ is the generator serial number, $\delta_i$ is the rotor angle,
$\omega_i$ is the rotor speed in rad/s, and $e'_{qi}$ and $e'_{di}$ are the transient voltage along $q$ and $d$ axes; $i_{qi}$ and $i_{di}$ are stator currents at $q$ and $d$ axes;
$T_{mi}$ is the mechanical torque, $T_{ei}$ is the electric air-gap torque, and $E_{fdi}$ is the internal field voltage; $\omega_0$ is the rated value of angular frequency, $H_i$ is the inertia constant, and $K_{Di}$ is the damping factor; $T'_{q0i}$ and $T'_{d0i}$ are the open-circuit time constants for $q$ and $d$ axes; $x_{qi}$ and $x_{di}$ are the synchronous reactance and $x'_{qi}$ and $x'_{di}$ are the transient reactance at the $q$ and $d$ axes.

Generators with classical model are described by the first two equations of (\ref{gen model}) and $e'_{qi}$ and $e'_{di}$ are kept unchanged. 

$T_{mi}$ and $E_{fdi}$ are considered as inputs and are assumed to be constant and known.
Let $\mathcal{G}_{P}$ denote the set of generators where PMUs are installed. For generator $i \in \mathcal{G}_{P}$, the terminal voltage phaosr $E_{ti}=e_{Ri}+je_{Ii}$ and the terminal current phasor $I_{ti}=i_{Ri}+ji_{Ii}$ can be measured, and are used as the outputs. 

The dynamic model (\ref{gen model}) can be rewritten in a general state space form in (\ref{n1}) and
the state vector $\boldsymbol{x}$, input vector $\boldsymbol{u}$, and output vector $\boldsymbol{y}$
can be written as
\begin{align}
\boldsymbol{x} &= \bmat{\boldsymbol{\delta}^\top \quad \boldsymbol{\omega}^\top \quad \boldsymbol{e'_{q}}^\top \quad \boldsymbol{e'_{d}}^\top}^\top  \\ 
\boldsymbol{u} &= \bmat{\boldsymbol{T_m}^\top \quad \boldsymbol{E_{fd}}^\top}^\top  \\
\boldsymbol{y} &= \bmat{\boldsymbol{e_R}^\top \quad \boldsymbol{e_I}^\top \quad \boldsymbol{i_R}^\top \quad \boldsymbol{i_I}^\top}^\top. 
\end{align}

The $i_{qi}$, $i_{di}$, and $T_{ei}$ in (\ref{gen model}) are actually functions of $\boldsymbol{x}$:
\begin{align*}
\it \Psi_{Ri}&=e'_{di}\sin\delta_i+e'_{qi}\cos\delta_i \\
\it \Psi_{Ii}&=e'_{qi}\sin\delta_i-e'_{di}\cos\delta_i \\
I_{ti}&=\overline{\boldsymbol{Y}}_i(\boldsymbol{\it \Psi}_{R}+j\boldsymbol{\it \Psi}_{I}) \\
i_{Ri}&= \operatorname{Re}(I_{ti}) \\
i_{Ii}&= \operatorname{Im}(I_{ti}) \\
i_{qi} &= \frac{S_B}{S_{Ni}}(i_{Ii}\sin\delta_i+i_{Ri}\cos\delta_i) \\
i_{di} &= \frac{S_B}{S_{Ni}}(i_{Ri}\sin\delta_i-i_{Ii}\cos\delta_i) \\
e_{qi}&=e'_{qi}-x'_{di}i_{di} \\
e_{di}&=e'_{di}+x'_{qi}i_{qi} \\
P_{ei} &= e_{qi}i_{qi}+e_{di}i_{di} \\
T_{ei} &= \frac{S_B}{S_{Ni}} P_{ei}
\end{align*}
where $\it \Psi_i=\Psi_{Ri}+j\Psi_{Ii}$ is the voltage source, $\boldsymbol{\it \Psi}_R$ and $\boldsymbol{\it \Psi}_I$ are column vectors of all generators' $\it \Psi_{Ri}$ and $\it \Psi_{Ii}$, $e_{qi}$ and $e_{di}$ are the terminal voltage at $q$ and $d$ axes, and $\overline{\boldsymbol{Y}}_i$ is the $i$th row of the admittance matrix of the reduced network $\boldsymbol{\overline{Y}}$ whose elements are constant if the difference between $x'_d$ and $x'_q$ is ignored (\cite{Wang1}), $P_{ei}$ is the electrical active output power, and $S_B$ and $S_{Ni}$ are the system base MVA and the base MVA for generator $i$, respectively.

The outputs $i_R$ and $i_I$ have been written as functions of $\boldsymbol{x}$.
Similarly, the outputs $e_{Ri}$ and $e_{Ii}$ can also be written as function of $\boldsymbol{x}$:
\begin{align*}
e_{Ri}&= e_{di}\sin\delta_i+e_{qi}\cos\delta_i \\
e_{Ii}&=e_{qi}\sin\delta_i-e_{di}\cos\delta_i.
\end{align*}

The continuous model in (\ref{n1}) can be discretized as 
\begin{equation}~\label{n2}
\left\{
\arraycolsep=1.4pt\def\arraystretch{1.4}
\begin{array}{ll}
\boldsymbol{x}_k = \boldsymbol{f}_d(\boldsymbol{x}_{k-1},\boldsymbol{u}_{k-1}) \\
\boldsymbol{y}_k = \boldsymbol{h}(\boldsymbol{x}_k,\boldsymbol{u}_k)
\end{array}
\right.
\end{equation}
where 
the state transition functions $\boldsymbol{f}_d$ can be obtained by the modified Euler method as
\begin{align}
\tilde{\boldsymbol{x}}_k &= \boldsymbol{x}_{k-1}+\boldsymbol{f}(\boldsymbol{x}_{k-1},\boldsymbol{u}_{k-1})\Delta t  \\
\tilde{\boldsymbol{f}} &= \frac{\boldsymbol{f}(\tilde{\boldsymbol{x}}_k,\boldsymbol{u}_k) + \boldsymbol{f}(\boldsymbol{x}_{k-1},\boldsymbol{u}_{k-1})}{2}\\
\boldsymbol{x}_k &= \boldsymbol{x}_{k-1}+\tilde{\boldsymbol{f}}\Delta t.
\end{align}

\section{Empirical Observability Gramian} \label{gramian}

Empirical observability gramian (\cite{Lall,Lall1,Hahn,Singh1}) provide a computable tool for
empirical analysis of the state-output behavior of nonlinear systems, 
which can be used to quantify the observability of the system state under some specific sensor placement.

The following sets are defined for empirical observability gramian:
\begin{align}
& T^n=\{\boldsymbol{T}_1,\cdots,\boldsymbol{T}_r;\;\;\;\boldsymbol{T}_l \in \mathbb{R}^{n\times n},\;\boldsymbol{T}_l^\top \boldsymbol{T}_l=\boldsymbol{I}_n,\;l=1,\ldots,r\} \nonumber \\
& M=\{c_1,\cdots,c_s;\;\;\;\;c_m \in \mathbb{R},\;c_m>0,\;m=1,\ldots,s\} \nonumber \\
& E^n=\{\boldsymbol{e}_1,\cdots,\boldsymbol{e}_n;\;\;\;\textrm{standard unit vectors in}\;\mathbb{R}^n\} \nonumber
\end{align}
where $T^n$ defines initial state perturbation directions, $r$ is the number of matrices for perturbation directions, 
$\boldsymbol{I}_n$ is an $n\times n$ identity matrix; $M$ defines perturbation sizes and $s$ is the number of different perturbation sizes 
for each direction, and $E^n$ defines the state to be perturbed and $n$ is the number of states.

For the nonlinear system described by (\ref{n1}), the empirical observability gramian can be defined as (\cite{Lall})
\begin{equation} \label{gd0}
\boldsymbol{W}_o^{con}=\sum_{l=1}^{r}\sum_{m=1}^{s}\frac{1}{rsc_m^2}\int_0^\infty \boldsymbol{T}_l \boldsymbol{\Psi}^{lm}(t)\boldsymbol{T}_l^\top dt
\end{equation}
where $\boldsymbol{\Psi}^{lm}(t)\in \mathbb{R}^{n\times n}$ is given by $\boldsymbol{\Psi}_{ij}^{lm}(t)=(\boldsymbol{y}^{ilm}(t)-\boldsymbol{y}^{ilm,0})^\top(\boldsymbol{y}^{jlm}(t)-\boldsymbol{y}^{jlm,0})$, $\boldsymbol{y}^{ilm}(t)$ is the output of the nonlinear system corresponding to the initial condition $\boldsymbol{x}(0)=c_m \boldsymbol{T}_l \boldsymbol{e}_i+\boldsymbol{x}_0$, and $\boldsymbol{y}^{ilm,0}$ refers to the output measurement corresponding to the unperturbed initial state $\boldsymbol{x}_0$, which is usually chosen as the steady state under typical power flow conditions.

The discrete form of empirical observability gramian can be defined as (\cite{Hahn})
\begin{equation} \label{gd}
\boldsymbol{W}_o=\sum_{l=1}^{r}\sum_{m=1}^{s}\frac{1}{rsc_m^2}\sum_{k=0}^K \boldsymbol{T}_l \boldsymbol{\Psi}^{lm}_k \boldsymbol{T}_l^\top \Delta t_k
\end{equation}
where $\boldsymbol{\Psi}^{lm}_k \in \mathbb{R}^{n\times n}$ is given by ${\boldsymbol{\Psi}^{lm}_k}_{ij}=(\boldsymbol{y}^{ilm}_k-\boldsymbol{y}^{ilm,0})^\top(\boldsymbol{y}^{jlm}_k-\boldsymbol{y}^{jlm,0})$, $\boldsymbol{y}^{ilm}_k$ is the output at time step $k$ corresponding to initial condition $\boldsymbol{x}(0)=c_m \boldsymbol{T}_l \boldsymbol{e}_i+\boldsymbol{x}_0$, $\boldsymbol{y}^{ilm,0}$ is the output for unperturbed initial state $\boldsymbol{x}_0$, $K$ is the number of points, and $\Delta t_k$ is the time interval.

\section{Adaptive Optimal PMU Placement} \label{adaptive}

Based on the empirical observability gramian in (\ref{gd}), the OPP for dynamic state estimation can be formulated as
\begin{align} \label{opt}
&\max\limits_{\boldsymbol{z}} F\big( \, \boldsymbol{W}_o(\boldsymbol{z})\big) \nonumber \\
\textrm{s.t.}\; & \;\;\; \sum_{i=1}^g  z_i = \bar{g} \notag \\
z_i & \in \{0,1\}, \;\; i=1,\ldots,g 
\end{align}
where $F$ is a function (measure) of the gramian $\boldsymbol{W}_o$, $\boldsymbol{z}$ is the vector of binary control variables determining PMU placement, and $g$ and $\bar{g}$ are, respectively, the number of generators and PMUs.

It has been shown that the empirical observability gramian for a system with $p$ outputs 
is the summation of the empirical gramians computed for each of the $p$ outputs individually (\cite{Singh1,Qi1}).
Obviously, the empirical observability gramian calculated from placing $\bar{g}$ PMUs individually
adds to be the identical gramian calculated from placing the $\bar{g}$ PMUs simultaneously.
Therefore, the optimization problem (\ref{opt}) can be rewritten as
\begin{align}\label{opt1}
&\max\limits_{\boldsymbol{z}} F\Big( \, \sum_{i=1}^g z_i \boldsymbol{W}_{o,i}(\boldsymbol{z})\Big) \nonumber \\
\;\;\textrm{s.t.}\; & \qquad \sum_{i=1}^g  z_i = \bar{g}  \nonumber \\
&\qquad z_i \in \{0,1\}, \;\; i=1,\ldots,g 
\end{align}
where $\boldsymbol{W}_{o,i}$ is the empirical observability gramian by only placing one PMU at generator $i$.

The objective function $F$ quantifies the observability of the system state and can be the determinant ($\det$) (\cite{Qi1,Qi4}), the trace ($\Tr$) (\cite{Singh1}), 
the smallest eigenvalue ($\sigma_{\min}$) (\cite{Kang,Kang1,Kang2,Sun1,Krener}), or the opposite of the condition number ($-\kappa$) (\cite{Krener}). By solving (\ref{opt1}) an OPP can be obtained for each $\bar{g}$ and each $F$.

Different measures of the observability gramian reflect various aspects of observability. The determinant and the trace measure the overall observability in all directions in noise space. However, the trace cannot tell the existence of a zero eigenvalue and an unobservable system may still have a large trace. The smallest eigenvalue defines the worst scenario of observability
while the condition number emphasizes the numerical stability in state estimation.

Although the determinant is a good measure that reflects the overall observability in all directions in noise space, 
it can be too small to be able to indicate the degree of observability, especially when the number of PMUs is small, in which case the observability for some states can be very weak 
and the determinant of the empirical observability gramian (the product of all eigenvalues) can be very small. 
When this is the case, the smallest eigenvalue or the condition number are better measures for quantifying the observability. 

Based on this, we propose an adaptive OPP method by making use of the advantages of different measures as:
\begin{displaymath}
\boldsymbol{z}^*(\bar{g}) = \left\{
\arraycolsep=1.4pt\def\arraystretch{1.4}
\begin{array}{ll}
\boldsymbol{z}_{\det}^*(\bar{g}) & \textrm{if $\det\big(\boldsymbol{z}_{\det}^*(\bar{g})\big) \ge \epsilon$} \\
\boldsymbol{z}_{-\kappa}^*(\bar{g}) & \textrm{if $\det\big(\boldsymbol{z}_{\det}^*(\bar{g})\big) < \epsilon$ and $R_{-\kappa} \ge R_{\sigma_{\min}}$} \\
\boldsymbol{z}_{\sigma_{\min}}^*(\bar{g}) & \textrm{if $\det\big(\boldsymbol{z}_{\det}^*(\bar{g})\big) < \epsilon$ and $R_{-\kappa} < R_{\sigma_{\min}}$}
\end{array} \right.
\end{displaymath}
where $\boldsymbol{z}^*(\bar{g})$ is the adaptive OPP for placing $\bar{g}$ PMUs, $\boldsymbol{z}_F^*(\bar{g})$ is the optimal placement of $\bar{g}$ PMUs by using $F\in \{\det,\sigma_{\min},-\kappa\}$ as the objective function in (\ref{opt1}), and 
\begin{align}
&R_{-\kappa}=\frac{\kappa\big(\boldsymbol{z}_{\det}^*(\bar{g})\big)}{\kappa\big(\boldsymbol{z}_{-\kappa}^*(\bar{g})\big)}\;, \\
&R_{\sigma_{\min}}=\frac{\sigma_{\min}\big(\boldsymbol{z}_{\sigma_{\min}}^*(\bar{g})\big)}{\sigma_{\min}\big(\boldsymbol{z}_{\det}^*(\bar{g})\big)}
\end{align}
are used to indicate the improvement of the condition number using $-\kappa$ objective function or the minimum eigenvalue using $\sigma_{\min}$ objective function compared with using $\det$ objective function, and  we write $F\Big(\sum_{i=1}^{g}z_{F,i}^*\boldsymbol{W}_{o,i}\big(\boldsymbol{z}^*_{F}(\bar{g})\big)\Big)$ in short as $F\big(\boldsymbol{z}^*_F(\bar{g})\big)$, which denotes the function value of $F\in \{\det,\sigma_{\min},\kappa \}$ for $\sum_{i=1}^{g}z_{F,i}^*(\bar{g})\boldsymbol{W}_{o,i}\big(\boldsymbol{z}^*_{F}(\bar{g})\big)$.
When the determinant is less than $\epsilon$, the observability is too weak and the determinant is too small to be a good measure of observability. 
In this paper $\epsilon$ is chosen as $1$.

In the proposed adaptive optimal PMU placement method, the optimization problem in (\ref{opt1}) has to be solved for three times. 
Compared with the method in \cite{Qi1}, this will increase the computational burden. 
However, since the PMU placement problem is used for planning and is solved offline, the calculation time is not a concern.

\section{Case Study} \label{case}

\begin{figure}[!t]
\centering
\includegraphics[width=3.2in]{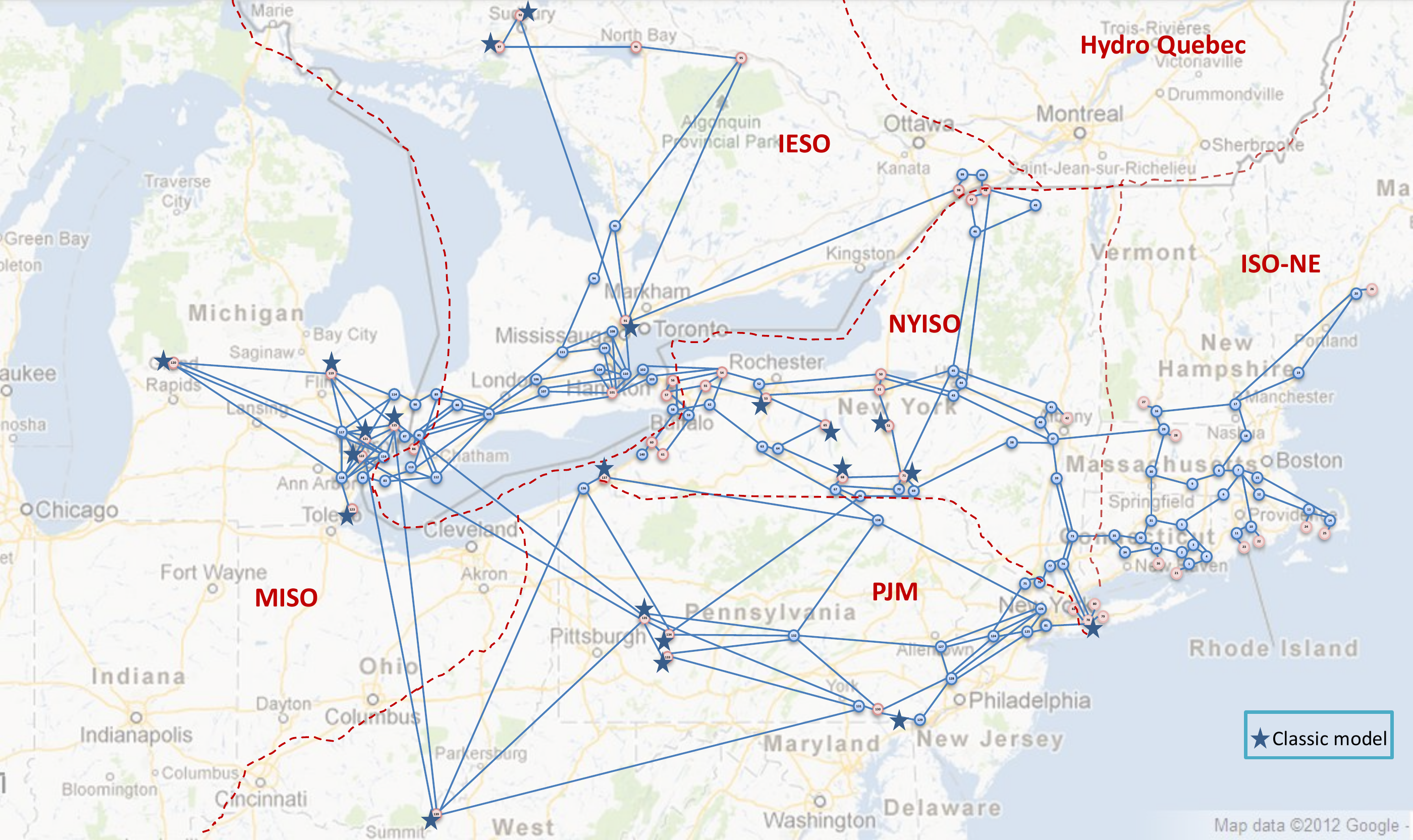}
\caption{Map of the NPCC 48-machine 140-bus system. The stars indicates generators with classical model.}
\label{npcc48}
\end{figure}

The proposed method is tested on the NPCC 48-machine 140-bus system (\cite{Chow}), which represents the northeast region of the Eastern Interconnection system. Among the 48 generators, 27 of them have fourth-order transient model and the others have second-order classical model. The map of this system is shown in Fig. \ref{npcc48}.

The empirical observability gramian calculation in (\ref{gd}) is implemented based on emgr (Empirical Gramian Framework) (\cite{Himpe}) on time interval $[0,t_f]$. In this paper $t_f$ is chosen as 5 seconds. 

The mixed-integer optimization problem in (\ref{opt1}) is solved by NOMAD solver (\cite{Digabel}),
which is a derivative-free global mixed integer nonlinear programming solver and is called by the OPTI toolbox (\cite{Currie}). 
NOMAD implements the Mesh Adaptive Direct Search (MADS) algorithm (\cite{Audet}),
a derivative-free direct search method with a rigorous convergence theory based on nonsmooth calculus (\cite{Clarke}), 
and aims for the best possible solution with a small number of evaluations. 

DSE is performed by using square-root unscented Kalman filter (\cite{Merwe,Qi1,Qi2}) for 50 times under OPP for different $F$ and random PMU placements to estimate the state trajectory on $[0,5s]$. For each case a three-phase fault is applied at the from bus of one of the 50 branches with highest line flows.
Measurements are voltage phasors and current phasors from PMUs at generator terminal buses with a sampling rate of 60 frames/s. All the other settings are the same as (\cite{Qi1}).

\subsection{Comparing PMU Placements Using Different Measures of the Empirical Observability Gramian} \label{result}

\begin{figure}
\begin{center}
\includegraphics[width=2.9in]{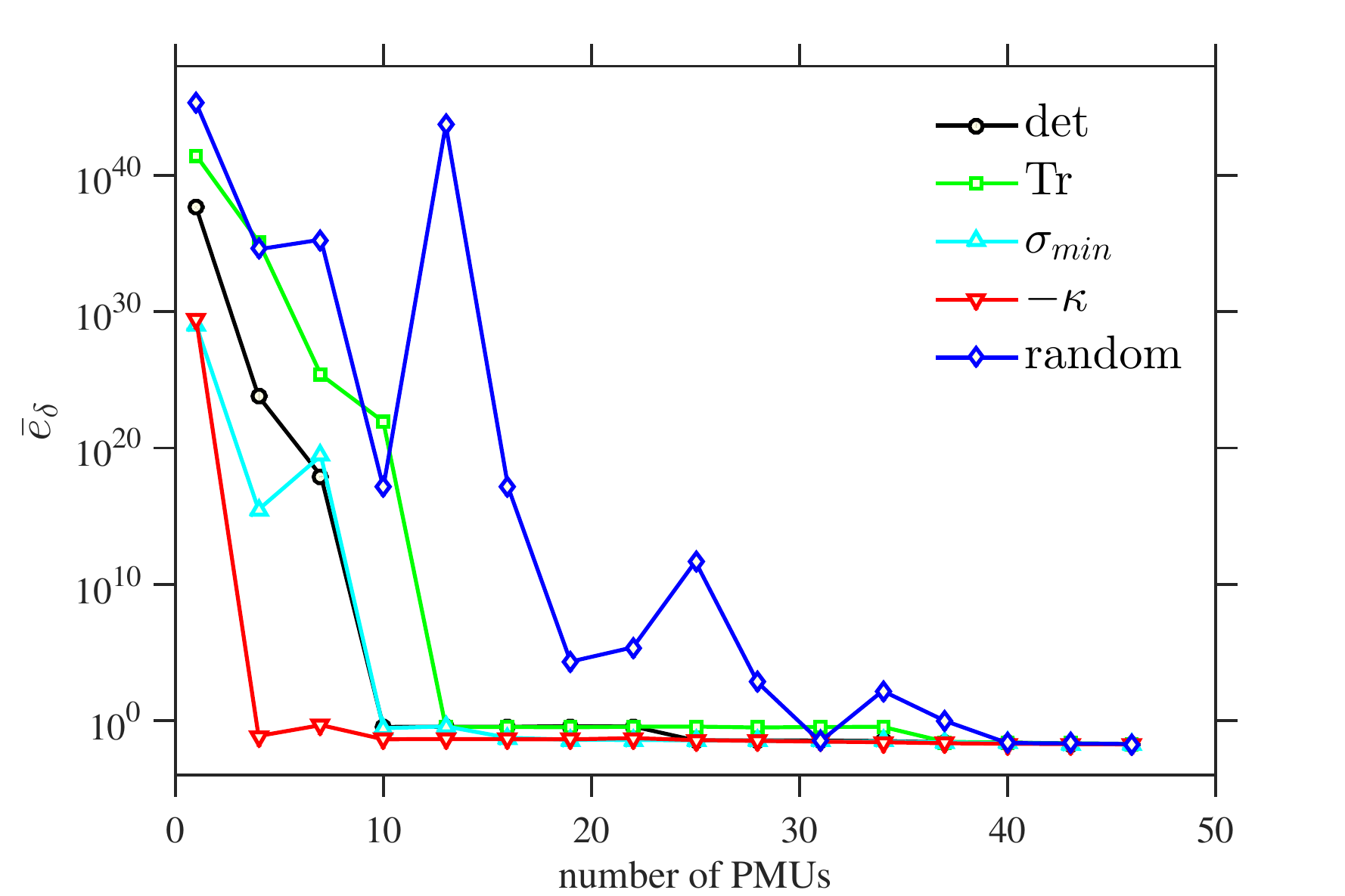}
\caption{Average error for different PMU placements} 
\label{error}
\end{center}
\end{figure}

In Fig. \ref{error} we show the average error of rotor angles
\begin{equation}
\bar{e}_\delta = \sqrt{\frac{\sum\limits_{i=1}^g \sum\limits_{t=1}^{T_s}\big(\delta_{i,t}^{\textrm{est}}-\delta_{i,t}^{\textrm{true}}\big)^2}{g\,T_s}}
\end{equation}
where $\delta_{i,t}^{\textrm{est}}$ and $\delta_{i,t}^{\textrm{true}}$ are the estimated and true rotor angles of the $i$th generator at time step $t$ and $T_s$ is the number of time steps. 

\begin{figure}[!t]
\centering
\includegraphics[width=2.9in]{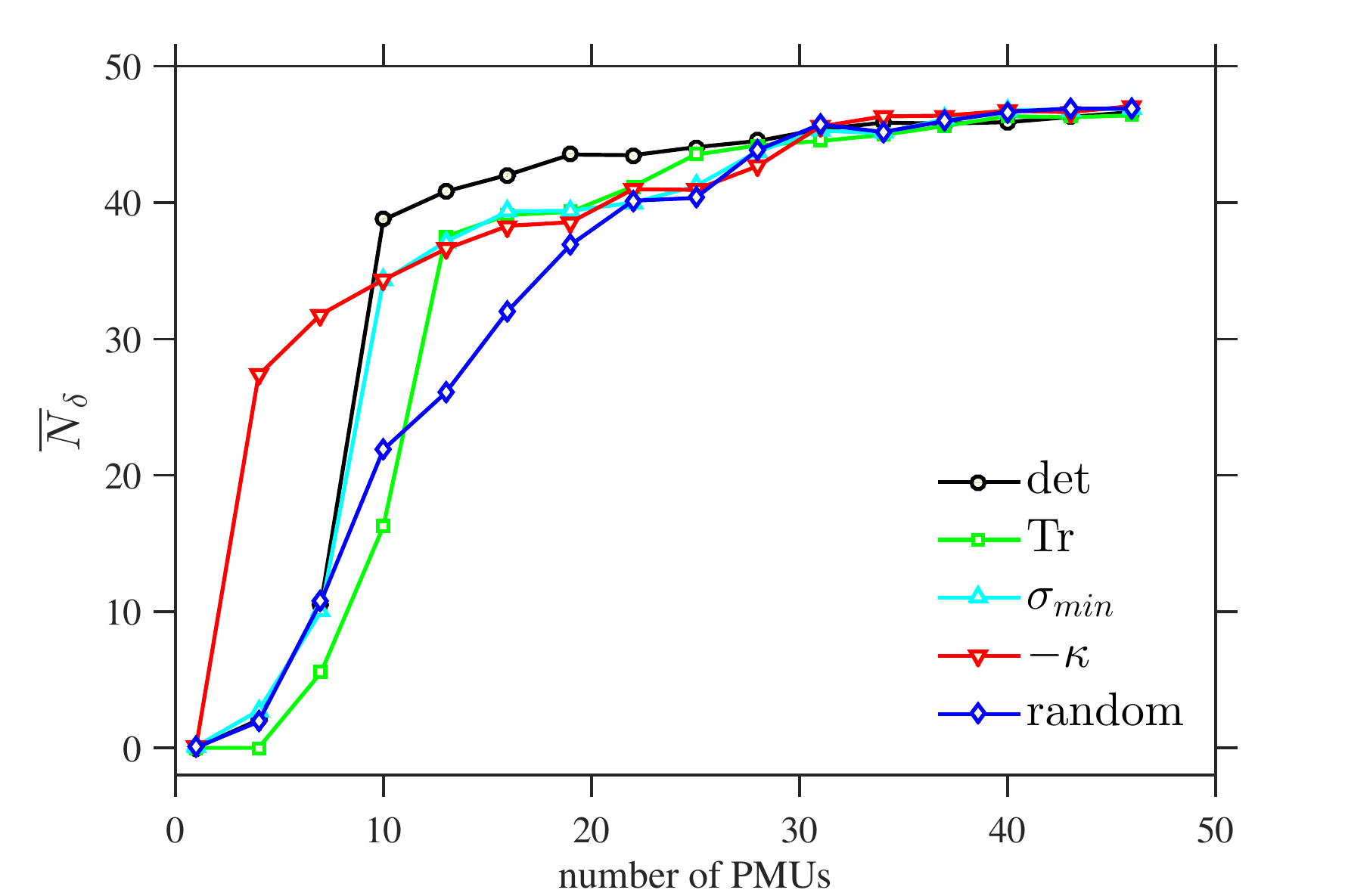}
\caption{Number of convergent angles under different PMU placements.}
\label{num}
\end{figure}

Fig. \ref{num} shows the average number of convergent angles, $\overline{N}_\delta$, for which the differences between the estimated and true values in the last one second are less than 2\% of the absolute value of true values.

From these figures it is seen that: 

\begin{enumerate}
\item Optimal placements are better than random placements in terms of guaranteeing better observability that can be indicated by smaller estimation error and greater number of convergent states.
\item $\Tr$ is not a good measure, mainly due to its being unable to tell the existence of a zero eigenvalue.
\item $\det$ is a good measure if the number of PMUs is not too small and the observability is not too weak; otherwise it has too small values to indicate observability.
\item $\sigma_{\min}$ and $-\kappa$ work better than $\det$ for small number of PMUs, in which case they have more reasonable values to indicate observability.
\end{enumerate}

\subsection{Results for Adaptive Optimal PMU Placement}

\begin{figure}[!t]
\centering
\includegraphics[width=3.0in]{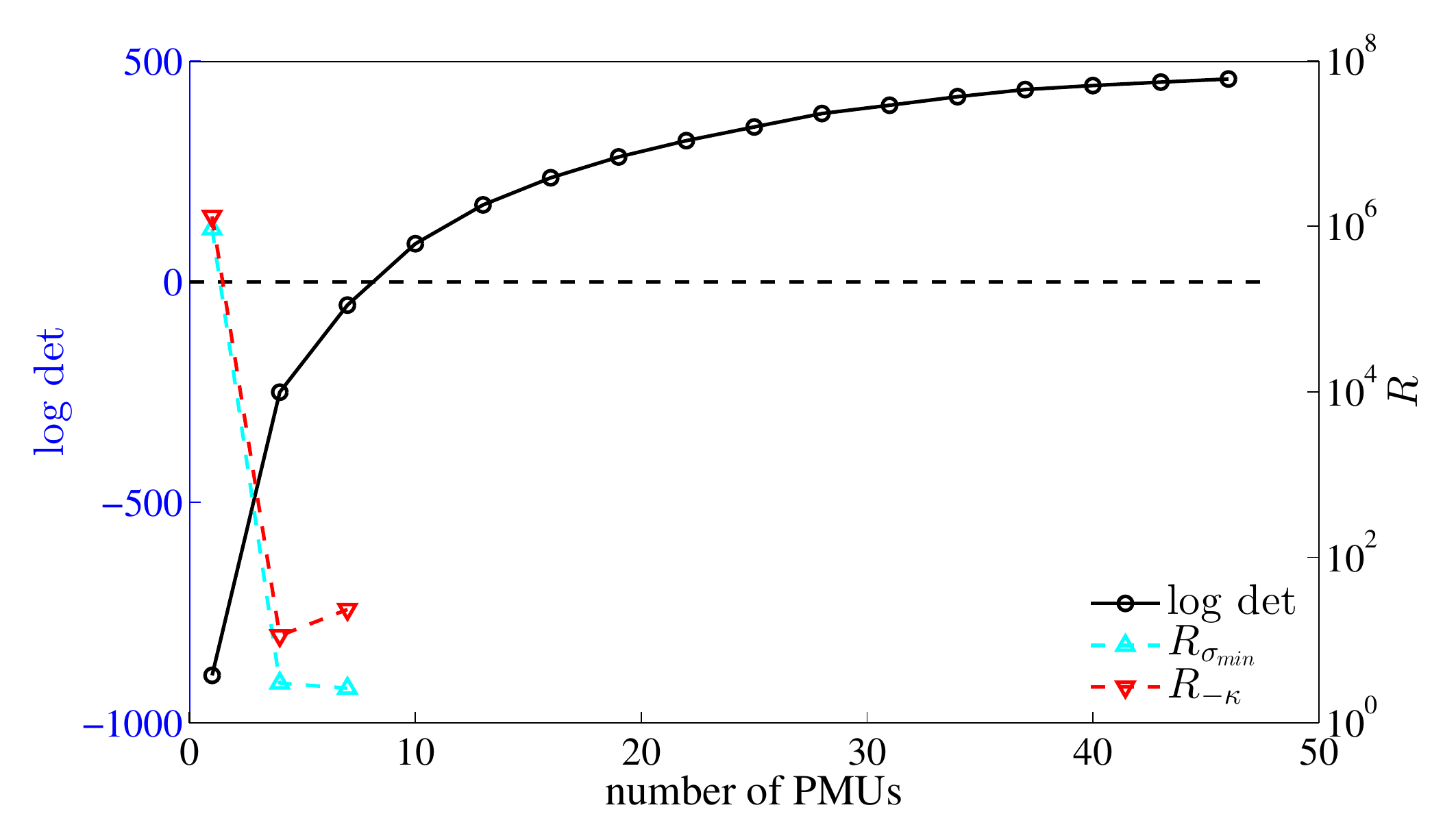}
\caption{Illustration of the adaptive OPP method. The dash line denotes that the maximum determinant of the empirical observability gramian equals 1 (or equivalently the corresponding $\log \det$ equals 0). }
\label{det}
\end{figure}

Fig. \ref{det} illustrates the proposed method. When the number of PMUs $\bar{g} \ge 10$, $\det(\boldsymbol{z}_{\det}^*) \ge 1$ ($\log \det(\boldsymbol{z}_{\det}^*) \ge 0$) and thus the determinant is big enough to be used as objective function. When $\bar{g} < 10$, $\sigma_{\min}$ or $-\kappa$ is considered. Since the improvement for condition number is higher ($R_{-\kappa} > R_{\sigma_{\min}}$), $-\kappa$ is chosen as the objective function.

\begin{figure}[!t]
\centering
\includegraphics[width=2.9in]{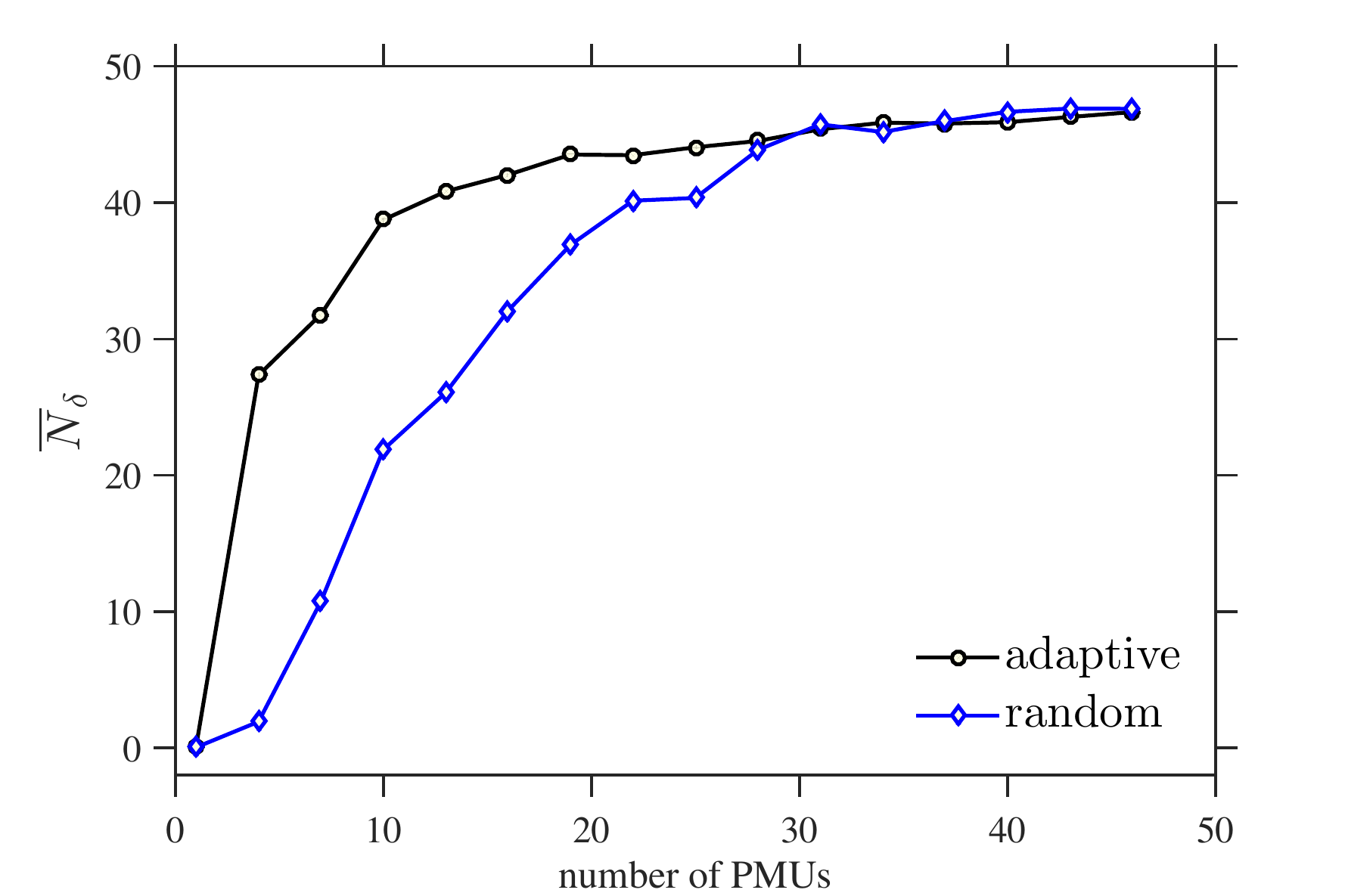}
\caption{Number of convergent angles for adaptive OPP compared with random OPP.}
\label{compare}
\end{figure}

In Fig. \ref{compare} we show the average numbers of convergent angles for adaptive OPP, 
which are much greater than those under random placements and are also better than those only using one measure as the objective function.

%
%

\section{Conclusion} \label{conclusion}

In this paper, different measures of the empirical observability gramian that can quantify the degree of observability of a power system under a specific PMU configuration are compared and an adaptive optimal PMU placement method for dynamic state estimation is proposed by automatically choosing proper measures as the objective function to guarantee best observability. 
For systems with a relatively small number of PMUs, the observabiity measured by condition number and smallest eigenvalue are more reliable, i.e. the numerical stability and the worst case scenario have higher impact on the observability. On the other hand, if the number of PMUs is not too small, the determinant is a better overall measure of observability.

The proposed method is tested on the NPCC 48-machine 140-bus system 
by performing dynamic state estimation with square-root unscented Kalman filter. 
Compared with only using one measure to quantify observability, the proposed method 
can guarantee smaller estimation errors and larger number of convergent states.

The discussion on different measures of the empirical observability gramian should be valid for the empirical controllability gramian. 
The idea based on which the adaptive optimal PMU placement method is proposed should also work for other placement problems related to either observability or controllability, 
such as the dynamic var sources placement discussed in \cite{Qi4}.

\end{document}